\documentclass[10pt]{article}
\usepackage{amsfonts}
\usepackage{amsmath}
\usepackage{mathrsfs}
\usepackage{mathrsfs,amscd,amssymb,amsthm,amsmath,bm,graphicx,psfrag,subfigure}

\makeatletter

\renewcommand{\@seccntformat}[1]{{\csname the#1\endcsname}{\normalsize .}\hspace{.5em}}
\makeatother

\def \[{\begin{equation}}
\def \]{\end{equation}}

\def \dist{{\rm dist}}
\newtheorem{thm}{Theorem}[section]

\newtheorem{defi}{Definition}
\newtheorem{claim}{Claim}
\newtheorem{lem}[thm]{Lemma}

\newtheorem{remark}{Remark}

\newenvironment{wst}
{\setlength{\leftmargini}{1.5\parindent}
 \begin{itemize}
 \setlength{\itemsep}{-1.1mm}}
{\end{itemize}}

\begin{document}

\setlength{\baselineskip}{15pt}
\begin{center}{\Large \bf On the spectral moments of trees with a given bipartition\footnote{Financially supported by the National Natural
Science Foundation of China (Grant Nos. 11071096, 11271149) and the Special Fund for Basic Scientific Research of Central Colleges (CCNU11A02015).}}

\vspace{4mm}

{\large Shuchao Li\footnote{E-mail: lscmath@mail.ccnu.edu.cn (S.C.
Li), 453919067@qq.com (J. Zhang)},\ \ Jiajia Zhang}\vspace{2mm}

Faculty of Mathematics and Statistics,  Central China Normal
University, Wuhan 430079, P.R. China
\end{center}

\noindent {\bf Abstract}: For two given positive integers $p$ and $q$ with $p\leqslant q$, we denote $\mathscr{T}_n^{p, q}=\{T: T$ is a tree of order $n$ with a $(p, q)$-bipartition\}. For a graph $G$ with $n$  vertices, let $A(G)$ be its adjacency matrix with eigenvalues $\lambda_1(G), \lambda_2(G), \ldots, \lambda_n(G)$ in non-increasing order. The number $S_k(G):=\sum_{i=1}^{n}\lambda_i^k(G)\,(k=0, 1, \ldots, n-1)$ is called the $k$th spectral moment of $G$. Let $S(G)=(S_0(G), S_1(G),\ldots, S_{n-1}(G))$ be the sequence of spectral moments of $G$. For two graphs $G_1$ and $G_2$, one has $G_1\prec_s G_2$ if for some $k\in \{1,2,\ldots,n-1\}$, $S_i(G_1)=S_i(G_2)\,(i=0,1,\ldots,k-1)$ and $S_k(G_1)<S_k(G_2)$ holds. In this paper, the last four trees, in the $S$-order, among $\mathscr{T}_n^{p, q}\,(4\leqslant p\leqslant q)$ are characterized.

\vspace{2mm} \noindent{\it Keywords}: Spectral moment; $S$-order; Tree; Bipartition

\vspace{2mm}


 {\setcounter{section}{0}
\section{\normalsize Introduction}

Up to isomorphism, all graphs considered here are finite, simple and connected. Undefined terminology and notation may be referred to \cite{B-Z10}.
Let $G=(V_G,E_G)$ be a simple undirected graph with $n$ vertices. $G-v$, $G-uv$ denote the graph obtained from $G$ by deleting vertex $v \in V_G$, or edge
$uv \in E_G$, respectively (this notation is naturally extended if more than one vertex or edge is deleted). Similarly,
$G+uv$ is obtained from $G$ by adding edge $uv \not\in E_G$. For $v\in V_G,$ let
$N_G(v)$ (or $N(v)$ for short) denote the set of all the adjacent vertices of $v$ in $G$ and $d(v)=|N_G(v)|$. A \textit{leaf} of $G$ is a vertex of degree one.

Let $A(G)$ be the adjacency matrix of a graph $G$, and let $\lambda_1(G),\lambda_2(G),\ldots,\lambda_n(G)$ be the eigenvalues of a graph $G$ in non-increasing order. The number $\sum_{i=1}^n\lambda_i^k(G) (k=0,1,\ldots,n-1)$ is called the $k$th \textit{spectral moment} of $G$, denoted by $S_k(G)$. We know from \cite{B-Z1} that $S_0=0,\, S_1=l, \, S_2=2m,\, S_3=6t$, where $n, \,l,\, m, \,t$ denote the number of vertices, the number of loops, the number of edges and the number of triangles, respectively. Let $S(G)=(S_0(G),S_1(G),\ldots,S_{n-1}(G))$ be the sequence of spectral moments of $G$. For two graphs $G_1$ and $G_2$, we shall write $G_1=_s G_2$ if $S_i(G_1)=S_i(G_2)$ for $i=0,1,\ldots,n-1$. Similarly, we have $G_1\prec_sG_2$ ($G_1$ comes before $G_2$ in the $S$-order) if for some $k\, (1\leqslant k\leqslant {n-1})$, we have  $S_i(G_1)=S_i(G_2)\, (i=0,1,\ldots,k-1)$ and $S_k(G_1)<S_k(G_2)$. we shall also write $G_1\preceq_sG_2$ if $G_1\prec_sG_2$ or $G_1=_sG_2$. $S$-order was used in producing graph catalogs (see \cite{D-M}), and for a more general setting of spectral moments one may be referred to \cite{C-G}.

Investigation on $S$-order of graphs attracts more and more researchers' attention. Cvetkovi\'{c} and Rowlinson \cite{D-I} studied the $S$-order of trees and unicyclic graphs and characterized the first and the last graphs, in the $S$-order, of all trees and all unicyclic graphs with given girth, respectively. Wu and Liu \cite{C-R-S1,6} determined the last $\lfloor \frac{d}{2}+1\rfloor$ and the last $\lfloor \frac{g}{2}+2\rfloor$ graphs, in the $S$-order, of all $n$-vertex trees with diameter $d\,(4\leqslant d\leqslant n-3)$ and  all $n$-vertex unicyclic graphs of girth $g\ (3\leqslant g\leqslant n-3),$ respectively. Wu and Fan \cite{D-F} determined the first and the last graphs in the $S$-order, of all unicyclic graphs and bicyclic graphs, respectively. Pan et al. \cite{C-R-S2} gave the first $\sum_{k=1}^{\lfloor\frac{n-1}{3}\rfloor}\left(\lfloor\frac{n-k-1}{2}\rfloor-k+1\right)$ graphs apart from a path, in the $S$-order, of all trees on $n$ vertices. Pan et al. \cite{G-H-L} determined the last and the second last quasi-tree, in the $S$-order, among the set \text{$\mathscr{L}(n,d_0)=\{G: G$ is a quasi-tree of order $n$ with $G-u_0$ being a tree and $d_{G}(u_0)=d_0\}$}, respectively.

In light of the information available on the spectral moments of graphs, it is natural to consider this problem for some other class of graphs, and the trees with a $(p, q)$-bipartition are a reasonable starting point for such a investigation. The $n$-vertex trees with a $(p, q)$-bipartition have been considered in \cite{101,14,Y-L-Z}, whereas to our best knowledge, the spectral moments of trees in $\mathscr{T}_n^{p, q}\,(4\leqslant p\leqslant q)$ were, so far, not considered. Here, we identified the last four trees, in the $S$-order, among $\mathscr{T}_n^{p, q}\,(4\leqslant p\leqslant q)$. For more recent results on the spectral moments of graphs, one may be referred to \cite{3,4,10,12,13}.

Given a connected bipartite graph $G$ with $n$ vertices, its vertex set can be partitioned into
two subsets $V_1$ and $V_2$, such that each edge joins a vertex in $V_1$ with a vertex in $V_2$. Suppose that $V_1$ has
$p$ vertices and $V_2$ has $q$ vertices, where $p+q = n$ with $p \leqslant q$. Then we say that $G$ has a $(p, q)$-\textit{bipartition}. For convenience, let $\mathscr{T}_n^{p, q}$ be the set of all $n$-vertex trees, each of which has a $(p, q)$-bipartition.

Throughout the text we denote by $P_n,\, K_{1, n-1}$ and $C_n$ the path, star and cycle on $n$ vertices, respectively. Let $U_n$ be a graph obtained from $C_{n-1}$ by attaching a leaf to one vertex of $C_{n-1}$, and let $E_4$ be a graph obtained by deleting an edge
from a complete graph $K_4$; $E_5$ be a graph obtained from two cycles $C_3$ and $C_3'$ of length 3 by identifying
one vertex of $C_3$ with one vertex of $C_3'$. The graphs $U_4,\, U_5,\, E_4$ and $E_5$ are depicted in Fig. 1.
\begin{figure}[h!]
\begin{center}
  \psfrag{a}{$U_4$}\psfrag{b}{$U_5$}\psfrag{c}{$E_4$}\psfrag{d}{$E_5$}
  \includegraphics[width=150mm] {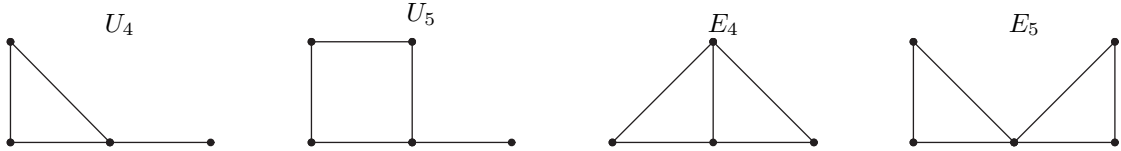}\\
  \caption{Four graphs $U_4,U_5,E_4$ and $E_5$.}
\end{center}
\end{figure}
Let $F$ be a graph. An $F$-\textit{subgraph} of $G$ is a subgraph of $G$ which is isomorphic to the graph $F$. Let $\phi_{G}(F)$ (or $\phi(F)$) be the number of all $F$-subgraphs of $G$. For a tree $T$ and
two vertices $v, u$ of $T$, the \textit{distance} $\dist_T(v, u)$ between $u$ and $v$ is the number of edges on the
unique path connecting them.

Further on we need the following lemmas.
\begin{lem}[\cite{D-I}]\label{lem1.1}
The $k$th spectral moment of $G$ is equal to the number of closed walks of length $k$.
\end{lem}

\begin{lem} \label{lem1.2}
For every graph $G$, we have
\begin{wst}
\item[{\rm (i)}] $S_4(G)=2\phi(P_2)+4\phi(P_3)+8\phi(C_4)$\, {\rm (see \cite{D-M}).}
\item[{\rm (ii)}] $S_5(G)=30\phi(C_3)+10\phi(U_4)+10\phi(C_5)$\, {\rm (see \cite{C-R-S1}).}
\item[{\rm (iii)}]
$S_6(G)=2\phi(P_2)+12\phi(P_3)+6\phi(P_4)+12\phi(K_{1,3})+12\phi(U_5)+36\phi(E_4)+24\phi(E_5)+24\phi(C_3)+48\phi(C_4)+12\phi(C_6)$\, {\rm (see \cite{C-R-S1}).}
\end{wst}
\end{lem}

Given a connected graph $G$, its line graph is denoted by $L(G)$. It is easy to see that the size of $L(G)$ is equal to the number of $P_3$ of $G$. By
[Exercise 1.5.10(a), 1], we have
\begin{lem}\label{lem1.3}
If $G$ is a simple connected graph, then
$\phi_G(P_3)=\sum_{v\in V_G}{d(v)\choose 2}.$
\end{lem}

\begin{defi}
Assume that $u, v, w$ are three distinct vertices of a tree $T$ satisfying $uv\in E_T,\, d(u)=1,\, d(w)\geqslant d(v)$ and $\dist_T(v,w)=2$. Let $T [v\rightarrow w;1]$ be the graph obtained from $T$ by deleting the edge $uv$ and adding the edge $uw$. In notation,
$$
\text{$T[v \rightarrow w;1]=T-uv+uw,$ }
$$
and we say $T[v \rightarrow w;1]$ is obtained from $T$ by
{Operation I}.
\end{defi}

\begin{remark}
If $T$ is in $\mathscr{T}_n^{p,q},$ by Definition 1, it is easy to see that $T[v\rightarrow w;1]$ is also in $\mathscr{T}_n^{p,q}.$
\end{remark}

\begin{lem}\label{lem1.4}
Let $T$ and $T[v \rightarrow w;1]$ be the trees defined as above. Then
$$
  T\prec_{s} T[v \rightarrow w;1].
$$
\end{lem}
\begin{proof}
By Lemma 1.1, $S_i(T)=S_i(T[v\rightarrow w;1])$ holds for $i=0,1,2,3$.
In view of Lemma 1.2(i), $\phi_{T}(P_2)=\phi_{T[v\rightarrow w;1]}(P_2)=n-1$, $\phi_{T}(C_{4})=\phi_{T[v\rightarrow w;1]}(C_{4})= 0$.
By Lemma 1.3, we have
$$
\phi_{T[v\rightarrow w;1]}(P_{3})- \phi_{T}(P_{3})={{d(w)+1} \choose 2}+{{d(v)-1} \choose 2}-{d(w)\choose 2}-{d(v)\choose 2}
                                                 =d(w)-d(v)+1>0.
$$
Hence, $S_{4}(T[v\rightarrow w;1])-S_{4}(T)=4(\phi_{T[v\rightarrow w;1]}(P_{3})-\phi_{T}(P_{3}))>0$, i.e.,
$T\prec_{s} T[v\rightarrow w;1]$.
\end{proof}

\begin{defi}
Let $uw$ be an edge of a tree $U$ with
$d(w)\geqslant 2$. $T$ is obtained from $U$ and the star $K_{1,k+1}$ ($k\geqslant2$) by
identifying $u$ with a pendant vertex of $K_{1,k+1}$ whose center is $v$.
Let $T[v\rightarrow w;2]$ be the graph obtained from $T$ by deleting all edges $vz$ and adding all
edges $wz$, where $z\in W=N_{T}(v)\backslash\{u\}.$ In notation,
$$
T[v\rightarrow w;2]=T-\{uz:z\in W\}+\{wz:z\in W\}
$$
and we say $T[v\rightarrow w;2]$ is obtained from $T$ by {Operation I\!\textrm{I}}.
Trees $T$ and $T[v\rightarrow w;2]$ are depicted in Fig. 2.
\end{defi}

\begin{remark}
If $T$ is in $\mathscr{T}_n^{p,q},$ by Definition 2, it is easy to see that $T[v\rightarrow w;2]$ is also in $\mathscr{T}_n^{p,q}.$
\end{remark}

\begin{figure}[h!]
\begin{center}
\psfrag{a}{$w$}\psfrag{b}{$v$} \psfrag{c}{$u$}\psfrag{d}{$v_1$}\psfrag{u}{$U$}\psfrag{f}{$v_k$} \psfrag{g}{$T$}
\psfrag{h}{$T[v\rightarrow w;2]$}
\psfrag{e}{$v_2$}
\includegraphics[width=90mm]{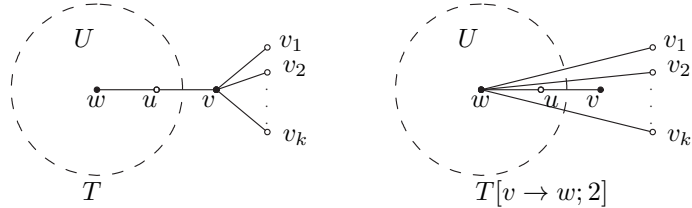}\\
\caption{$T\Rightarrow T[v\rightarrow w;2]$ by Operation I\!\textrm{I}.}
\end{center}
\end{figure}

\begin{lem}\label{thm1.5}
Let $T$ and $T[v\rightarrow w;2]$ be the trees described as above, one has $T\prec_s T[v\rightarrow w;2]$.
\end{lem}
\begin{proof}
By Lemma 1.1, $S_i(T)=S_i(T[v\rightarrow w;2])$ holds for $i=0,1,2,3$. In view of Lemma 1.2(i),
$\phi_{T}(P_2)=\phi_{T[v\rightarrow w;2]}(P_2)=n-1$ and $\phi_{T}(C_4)=\phi_{T[v\rightarrow w;2]}(C_4)= 0.$ By Lemma 1.3,
$$
\phi_{T[v\rightarrow w;2]}(P_3)- \phi_{T}(P_3)={{d(w)+k} \choose 2}-{d(w)\choose 2}-{{k+1}\choose 2}
                                                 =k(d(w)-1)>0.
$$
Hence, we have $S_4(T[v\rightarrow w;2])-S_4(T)=4(\phi_{T[v\rightarrow w;2]}(P_3)-\phi_{T}(P_3))>0$, i.e.,
$T\prec_s T[v\rightarrow w;2]$.
\end{proof}
\begin{figure}[h!]
\begin{center}
  \psfrag{a}{$v_0$}\psfrag{b}{$v_1$}\psfrag{c}{$v_2$}\psfrag{d}{$v_{i-1}$}\psfrag{e}{$v_i$}\psfrag{f}{$v_{i+1}$}\psfrag{x}{$v_{l-2}$}\psfrag{g}{$v_{l-1}$}
  \psfrag{h}{$v_l$} \psfrag{i}{$T_2$}\psfrag{j}{$T_{i-1}$}\psfrag{k}{$T_i$}\psfrag{l}{$T_{i+1}$}\psfrag{y}{$T_{l-2}$}\psfrag{r}{$T$}\psfrag{s}{$\hat{T}$}
  \psfrag{m}{$T'_2$}\psfrag{n}{$T'_{i-1}$}\psfrag{o}{$T'_{i+1}$}\psfrag{p}{$T'_i$}\psfrag{z}{$T'_{l-2}$}\psfrag{q}{$v_i$}\psfrag{t}{$T'_i$}
  \includegraphics[width=40mm] {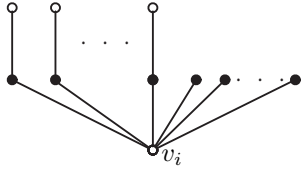}\\
  \caption{Tree $T_i'$.}
\end{center}
\end{figure}

For any $T\in \mathscr{T}_n^{p,q}$, $T$ contains a longest path, say $P=v_0v_1v_2\ldots v_l\,(l\geqslant 3)$.  Thus, $T$ can be seen as the tree  obtained from $P$ by attaching a tree $T_i$ to the vertex $v_i$ on $P$ if $d_T(v_i)>2$, $i\in \{2,3, \ldots, l-2\}$. Repeated applying {Operations I and I\!I} to $T$ yields an $n$-vertex tree, say $\hat{T}$, such that $P$ is also one of the longest paths of $\hat{T}$ and $d_{\hat{T}}(u)\leqslant 2$ for $u\in V_{\hat{T}}\setminus\{V_P\}.$ That is to say, $\hat{T}$ is obtained from $P$ by attaching a tree $T'_i$ to $v_i$ on $P$ if $d_{\hat{T}}(v_i)>2,\, i\in \{2,3, \ldots, l-2\},$ where $T_i'$ is depicted in Fig. 3. For convenience, let $\hat{\mathscr{T}}_n^{p,q}$ be the set of all trees $\hat{T}$ defined as above. By Remarks 1 and 2, $\hat{\mathscr{T}}_n^{p,q}$ is a proper subset of $\mathscr{T}_n^{p,q}.$

The following result follows immediately from Lemmas 1.4 and 1.5.
\begin{lem}\label{lem1.6}
For any $T\in \mathscr{T}_n^{p,q}\setminus \hat {\mathscr{T}}_n^{p,q}$, there exists a $\hat{T}\in \hat {\mathscr{T}}_n^{p,q}$ such that $T\prec_s \hat{T}$.
\end{lem}

Let $D_G^i(u)=\{v\in V_G: \dist_G(u,v)=i\}$,\, $i=2,3.$ Assume that $P=v_0v_1v_2\ldots v_l\ (l\geqslant  3)$ is one of the longest paths in a tree $T\in \hat{\mathscr{T}}_n^{p,q}$ and let $T'=T-\{v: \, v\in D_T^2(v_2)\bigcap D_T^3(v_3)\}+\{v_3v:\, v\in D_T^2(v_2)\bigcap D_T^3(v_3)\},$ where $T$ and $T'$ are depicted in Fig. 4. It is easy to see that $T'$ is obtained from $T$ by applying {Operation I}\, $a$ times and {Operation I\!I} once, by Lemmas 1.4 and 1.5, the following result holds immediately.
\begin{figure}[h!]
\begin{center}
  \psfrag{a}{$v_0$}\psfrag{b}{$v_1$}\psfrag{c}{$v_2$}\psfrag{d}{$v_3$}\psfrag{e}{$v_{l-2}$}\psfrag{f}{$v_{l-1}$}\psfrag{g}{$v_l$}\psfrag{h}{$T$}
  \psfrag{i}{$T'$} \psfrag{j}{$v_4$} \psfrag{1}{$a$} \psfrag{2}{$r$} \psfrag{3}{$s$} \psfrag{4}{$t$} \psfrag{5}{$t+a+r+1$}
  \includegraphics[width=120mm] {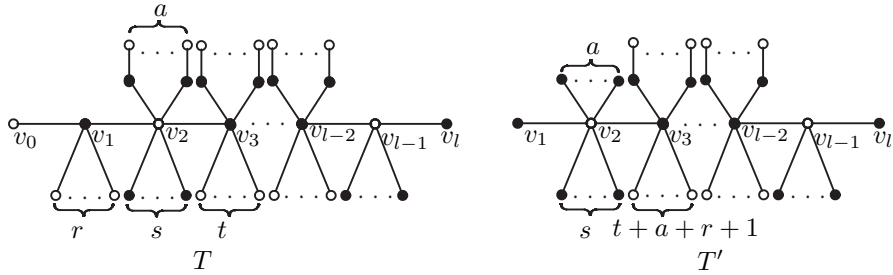}\\
  \caption{Trees $T$ and $T'$ with some vertices labeled.}
\end{center}
\end{figure}
\begin{lem}\label{lem1.7}
Let $T$ and $T'$ be the trees defined as above, one has $T\prec_s T'.$
\end{lem}

\section{\normalsize The last four trees in the $S$-order among $\mathscr{T}_n^{p,q}$ }\setcounter{equation}{0}
In this section, we determine the last four trees, in the $S$-order, among the set $\mathscr{T}_n^{p, q}\,(4\leqslant p\leqslant q)$.

For convenience, let
$B_{p,q}^{k,l},\, D_{p,q}^{k,l}\,(k,\,l\geqslant0)$ be the trees as depicted in Fig. 5, where the degree of $u$ is no less than that of $v$. In particular, $B_{p,q}^{0,0}\cong D_{p,q}^{0,0}.$
\begin{figure}
\begin{center}
  \psfrag{1}{$l$}\psfrag{2}{$k$}\psfrag{3}{$B_{p,q}^{k,l}$}\psfrag{4}{$D_{p,q}^{k,l}$}\psfrag{5}{$l$}\psfrag{6}{$k$}
  \psfrag{7}{$u$}\psfrag{8}{$v$}\psfrag{a}{$x$}\psfrag{b}{$y$}
  \includegraphics[width=100mm]{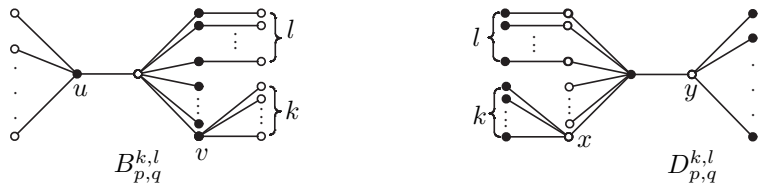}\\
  \caption{Trees $B_{p,q}^{k,l}$ and $D_{p,q}^{k,l}$ each of which contains $p$ white points and $q$ black points.}
\end{center}
\end{figure}
\begin{thm}\label{lem2.1}
Let $T$ be in $\mathscr{T}_n^{p,q}$, one has $T\preceq_s B_{p,q}^{0,0}$ with equality if and only if $T\cong B_{p,q}^{0,0}$.
\end{thm}
\begin{proof}
On the one hand, by Lemma 1.6, the last tree, in the $S$-order, among $\mathscr{T}_n^{p,q}$ must be in $\hat{\mathscr{T}}_n^{p,q}$. On the
other hand, for all $T\in \hat{\mathscr{T}}_n^{p,q},$ it is easy to see the tree $T'$ obtained from $T$, as depicted in Fig. 4, is also in
$\hat{\mathscr{T}}_n^{p,q}.$ Furthermore, the diameter of $T$ is $l;$ whereas the diameter of $T'$ is $l-1$; see Fig. 4. Applying this procedure (based on Operations I and I\!I) to $T'$ yields a tree, say $T''$. We have $T''\in \hat{\mathscr{T}}_n^{p,q}$, $T'\prec_s T''$ and the diameter of $T''$ is $l-2$. Repeat this procedure we finally obtain the unique graph $B_{p,q}^{0,0}$ of diameter 3 such that for any $T\in \hat{\mathscr{T}}_n^{p,q}\setminus \{B_{p,q}^{0,0}\}$ one has $T\prec_s B_{p,q}^{0,0}$ (based on
Lemma 1.7).
\end{proof}
\begin{thm}\label{lem2.2}
For any $T\in \mathscr{T}_n^{p,q}\setminus \{B_{p,q}^{0,0}\}$ with $4\leqslant p\leqslant q,$ one has $T\preceq_s B_{p,q}^{0,1}$ with equality if and only if $T\cong B_{p,q}^{0,1}$.
\end{thm}
\begin{proof}
For any $T \in \mathscr{T}_n^{p,q}$ such that $T\not\cong  B_{p,q}^{0,0}$, from the proof of Theorem 2.1, it is easy to see that $T$ can be transformed into $B_{p,q}^{0,0}$ by carrying the Operations I and I\!I repeatedly. Let $\mathscr{A}_1$ denote the set of all trees in $\mathscr{T}_n^{p,q}$ which can
be transformed into $B_{p,q}^{0,0}$ by carrying Operation I once, and let $\mathscr{A}_2$ denote the set of all trees in $\mathscr{T}_n^{p,q}$ which can be
transformed into $B_{p,q}^{0,0}$ by carrying Operation I\!I once. It follows from Lemmas 1.4 and 1.5 that the second
last tree, in the $S$-order, among $\mathscr{T}_n^{p,q}$ must be in $\mathscr{A}_1\cup \mathscr{A}_2$.

By definitions of $\mathscr{A}_1$ and $\mathscr{A}_2$, it is routine to check that $\mathscr{A}_1=\{B_{p,q}^{0,1},\,D_{p,q}^{0,1}\}$ (In particular, if $p=q$ then $B_{p,q}^{0,1}\cong D_{p,q}^{0,1}$; hence $\mathscr{A}_1=\{B_{p,q}^{0,1}\}$ for $p=q$), $\mathscr{A}_2=\{B_{p,q}^{k,0}: 2\leqslant k\leqslant \lfloor \frac{p-1}{2}\rfloor\}\cup \{D_{p,q}^{k,0}: 2\leqslant k\leqslant \lfloor \frac{q-1}{2}\rfloor\}$. Note that $B_{p,q}^{0,1}$ can be obtained from $B_{p,q}^{k,0}$ by using {Operation I}\, $(k-1)$ times, by Lemma 1.4, we have $B_{p,q}^{k,0}\prec_s B_{p,q}^{0,1}$ for $2\leqslant k\leqslant \lfloor \frac{p-1}{2}\rfloor$. Similarly, we have $D_{p,q}^{k,0}\prec_s D_{p,q}^{0,1}$ with $2\leqslant k\leqslant \lfloor \frac{q-1}{2}\rfloor$.

Hence, if $p=q$ then $B_{p,q}^{0,1}$ is just the second last tree, in the $S$-order, among $\mathscr{T}_n^{p,q}$ for $p\geqslant 4$. So in what follows we consider $p<q.$

In order to complete the proof, it suffices to compare $B_{p,q}^{0,1}$ with $D_{p,q}^{0,1}$.
By Lemma 1.1, we have $S_i(B_{p,q}^{0,1})=S_i(D_{p,q}^{0,1})$ for $i=0,1,2,3$. In view of Lemma 1.2(i), $\phi_{B_{p,q}^{0,1}}(P_2)=\phi_{D_{p,q}^{0,1}}(P_2)=n-1$ and $\phi_{B_{p,q}^{0,1}}(C_4)=\phi_{D_{p,q}^{0,1}}(C_4)=0$. In view of Lemma 1.3, we have
$$
\phi_{B_{p,q}^{0,1}}(P_3)-\phi_{D_{p,q}^{0,1}}(P_3)={{p-1} \choose 2}+{q \choose 2}+1-{p\choose 2}-{{q-1}\choose 2}-1
                                         =q-p>0.
$$
Hence, $S_4(B_{p,q}^{0,1})-S_4(D_{p,q}^{0,1})=4(\phi_{B_1}(P_3)-\phi_{D_{p,q}^{0,1}}(P_3))>0$, i.e., $D_{p,q}^{0,1}\prec_s B_{p,q}^{0,1}$.

This completes the proof.
\end{proof}

\begin{figure}[h!]
\begin{center}
  \psfrag{1}{$k$}\psfrag{2}{$C_{p,q}^k$}\psfrag{3}{$E_{p,q}^k$}
  \includegraphics[width=120mm]{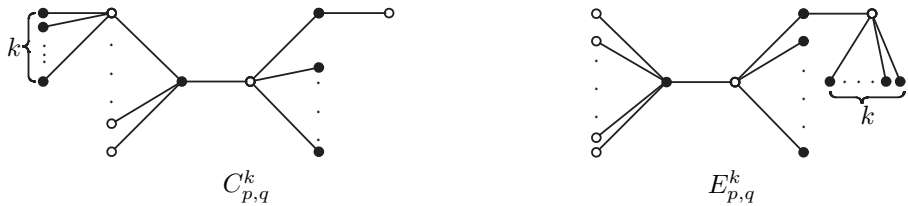}\\
  \caption{Trees $C_{p,q}^k$ and $E_{p,q}^k$ each of which contains $p$ white points and $q$ black points}
\end{center}
\end{figure}

For convenience, let
$C_{p,q}^k,\, E_{p,q}^k\,(1\leqslant k\leqslant q-2)$ be the trees as depicted in Fig. 6, it is easy to see that $C_{p,q}^k,\, E_{p,q}^k \in \mathscr{T}_n^{p,q}$.
\begin{thm}\label{lem2.4}
Let $p,q$ be positive integers with $4\leqslant p\leqslant q.$
\begin{wst}
  \item[{\rm (i)}] For any $T\in\mathscr{T}_n^{p,q}\setminus \{B_{p,q}^{0,0},\, B_{p,q}^{0,1}\}$ with $p=q$, we have $T\preceq_s B_{p,q}^{2,0}$  with equality if and only if $T\cong B_{p,q}^{2,0}$.
  \item[{\rm (ii)}] For any $T\in\mathscr{T}_n^{p,q}\setminus \{B_{p,q}^{0,0},\, B_{p,q}^{0,1}\}$ with $p<q$, if $p> \frac{q+4}{2}$, then  we have $T\preceq_s D_{p,q}^{0,1}$ with equality if and only if $T\cong D_{p,q}^{0,1};$ if $p\leqslant \frac{q+4}{2}$, then  we have $T\preceq_s B_{p,q}^{2,0}$  with equality if and only if $T\cong B_{p,q}^{2,0}$.
\end{wst}
\end{thm}
\begin{proof}
For any $T \in \mathscr{T}_n^{p,q}$ such that $T\not\cong  B_{p,q}^{0,0},\, B_{p,q}^{0,1}$, by a similar discussion as above, $T$ can be transformed into $B_{p,q}^{0,0}$ (resp. $B_{p,q}^{0,1}$) by carrying Operations I and I\!I repeatedly. Let $\mathscr{B}_1$ denote the set of all trees in $\mathscr{T}_n^{p,q}$ which can be transformed into $B_{p,q}^{0,1}$ by carrying Operation I once, and let $\mathscr{B}_2$ denote the set of all trees in $\mathscr{T}_n^{p,q}$ which can be transformed into $B_{p,q}^{0,1}$ by carrying Operation I\!I once. It follows from Lemmas 1.4 and 1.5, if $p<q$ then the third
last tree, in the $S$-order, among $\mathscr{T}_n^{p,q}$ must be in $\{D_{p,q}^{0,1}\}\cup\mathscr{A}_2\cup \mathscr{B}_1\cup \mathscr{B}_2$, where $\mathscr{A}_2$ is defined in the proof of Theorem 2.2. Note that if $p=q$, then $B_{p,q}^{0,1}\cong D_{p,q}^{0,1}$. Hence, the third
last tree, in the $S$-order, among $\mathscr{T}_n^{p,q}$ with $p=q$ must be in $\mathscr{A}_2\cup \mathscr{B}_1\cup \mathscr{B}_2.$

By the definition of $\mathscr{B}_1$ and $\mathscr{B}_2$, it's routine to check that $\mathscr{B}_1=\{B_{p,q}^{2,0},\,B_{p,q}^{0,2},\,C_{p,q}^1,\,E_{p,q}^1\}$ and $\mathscr{B}_2=\{C_{p,q}^k: 2\leqslant k\leqslant q-2\}\cup \{E_{p,q}^k: 2\leqslant k \leqslant q-2\}$. We first show the following two claims.
\begin{claim}
The last tree, in the $S$-order, among $\mathscr{A}_2$ is $B_{p,q}^{2,0}$.
\end{claim}
{\noindent\bf Proof of Claim 1}\ \
In graph $B_{p,q}^{k,0}$, we assume, without loss of generality, that $d(u)\geqslant  d(v)$ (see Fig. 5), hence $2\leqslant k\leqslant \lfloor\frac{p-1}{2}\rfloor$. Furthermore, we obtain (based on Lemma 1.4) that
\[
  B_{p,q}^{\lfloor\frac{p-1}{2}\rfloor,0}\prec_s B_{p,q}^{\lfloor\frac{p-1}{2}\rfloor-1,0}\prec_s \cdots \prec_s B_{p,q}^{k,0}\prec_s \cdots \prec_s B_{p,q}^{3,0} \prec_s B_{p,q}^{2,0}.
\]

Similarly, we obtain
\[
  D_{p,q}^{\lfloor\frac{q-1}{2}\rfloor,0}\prec_s D_{p,q}^{\lfloor\frac{q-1}{2}\rfloor-1,0}\prec_s \cdots \prec_s D_{p,q}^{k,0}\prec_s \cdots \prec_s D_{p,q}^{3,0} \prec_s D_{p,q}^{2,0}.
\]

Note that if $p=q$, it is easy to see that $B_{p,q}^{2,0}\cong D_{p,q}^{2,0},$ hence Claim 1 holds immediately. In what follows, we consider $p<q.$

In view of (2.1) and (2.2), it suffices to compare $B_{p,q}^{2,0}$ with that of $D_{p,q}^{2,0}$. In fact,
by Lemma 1.1 one has $S_i(B_{p,q}^{2,0})=S_i(D_{p,q}^{2,0})$ for $i=0,1,2,3$.
In view of Lemma 1.2(i), $\phi_{B_{p,q}^{2,0}}(P_2)=\phi_{D_{p,q}^{2,0}}(P_2)=n-1,\, \phi_{B_{p,q}^{2,0}}(C_4)=\phi_{D_{p,q}^{2,0}}(C_4)=0$ and by Lemma 1.3,
$$
\phi_{B_{p,q}^{2,0}}(P_3)-\phi_{D_{p,q}^{2,0}}(P_3)={{p-2} \choose 2}+{q \choose 2}-{p\choose 2}-{{q-2}\choose 2}
                                         =2(q-p)>0.
$$
Hence, we have $S_4(B_{p,q}^{2,0})-S_4(D_{p,q}^{2,0})>0$, i.e., $D_{p,q}^{2,0}\prec_s B_{p,q}^{2,0}$.

This completes the proof. \qed
\begin{claim}
The last tree, in the $S$-order, among $\mathscr{B}_1\cup \mathscr{B}_2$ is $B_{p,q}^{2,0}$.
\end{claim}
{\noindent\bf Proof of Claim 2}\ \
Note that $C_{p,q}^1$ can be obtained from $C_{p,q}^k$ by using Operation I\ $(k-1)$ times, hence by Lemma 1.4 $C_{p,q}^k\prec_s C_{p,q}^1$ for $k\geqslant  2$. Similarly, we have $E_{p,q}^k\prec_s E_{p,q}^1$ for $k\geqslant  2$. So the last tree, in the $S$-order, among $\mathscr{B}_1\cup \mathscr{B}_2$ must be in $\mathscr{B}_1$.

Note that $C_{p,q}^1$ and $E_{p,q}^1$ have the same degree sequence, by Lemma 1.3 we have
\[\label{eq:2.3}
 \phi_{E_{p,q}^1}(P_{3})=\phi_{C_{p,q}^1}(P_{3})
\]holds.
By Lemma 1.1, $S_i(B_{p,q}^{0,2})=S_i(C_{p,q}^1)=S_i(E_{p,q}^1)$ holds for $i=0,1,2,3$. In view of Lemma 1.2(i), it is routine to check that $\phi_{C_{p,q}^1}(P_2)=\phi_{E_{p,q}^1}(P_2)=\phi_{B_{p,q}^{0,2}}(P_2)=n-1$, $\phi_{C_{p,q}^1}(C_{4})=\phi_{E_{p,q}^1}(C_4)=\phi_{B_{p,q}^{0,2}}(C_{4})=0$. By Lemma 1.3, one has
$$
\phi_{C_{p,q}^1}(P_{3})-\phi_{B_{p,q}^{0,2}}(P_{3})=\left({{p-1}\choose 2}+{{q-1}\choose 2}+2\right)-\left({{p-2} \choose 2}+{q \choose 2}+2\right)                                     =-(q-p+1)<0.
$$
In view of (\ref{eq:2.3}) we have
$\phi_{C_{p,q}^1}(P_{3})-\phi_{B_{p,q}^{0,2}}(P_{3})<0$.
Hence, by Lemma 1.2(i), we have $S_{4}(C_{p,q}^1)<S_{4}(B_{p,q}^{0,2})$ and $S_{4}(E_{p,q}^1)<S_{4}(B_{p,q}^{0,2})$, i.e., $C_{p,q}^1\prec_{s}B_{p,q}^{0,2}$
and $E_{p,q}^1\prec_{s}B_{p,q}^{0,2}$.

On the other hand, $B_{p,q}^{0,2}$ can be transformed into $B_{p,q}^{2,0}$ by carrying {Operation I} once, by Lemma 1.4 we have $B_{p,q}^{0,2}\prec_{s}B_{p,q}^{2,0}$. That is to say,  $B_{p,q}^{2,0}$ is the last tree, in the $S$-order, among $\mathscr{B}_1\cup \mathscr{B}_2$.\qed

\medskip

If $p=q$, by Claims 1 and 2, we obtain that $B_{p,q}^{2,0}$ is just the last tree, in the $S$-order, among $\mathscr{T}_n^{p,q}\setminus \{B_{p,q}^{0,0},\, B_{p,q}^{0,1}\}.$ This completes the proof of (i).\medskip

Now we consider $p<q$ in what follows. According to Claims 1 and 2, it suffices to compare $B_{p,q}^{2,0}$ with $D_{p,q}^{0,1}$ in this case.

By Lemma 1.1, $S_i(D_{p,q}^{0,1})=S_i(B_{p,q}^{2,0})$ holds for $i=0,1,2,3$. In view of Lemma 1.2(i), it is routine to check that $\phi_{D_{p,q}^{0,1}}(P_2)=\phi_{B_{p,q}^{2,0}}(P_2)=n-1,\, \phi_{D_{p,q}^{0,1}}(C_4)=\phi_{B_{p,q}^{2,0}}(C_4)=0$. Furthermore, by Lemma 1.3, we have
\[\label{eq:2.4}
\phi_{D_{p,q}^{0,1}}(P_3)-\phi_{B_{p,q}^{2,0}}(P_3)=\left({p\choose 2}+{{q-1}\choose 2}+1\right)-\left({{p-2}\choose 2}+{q\choose 2}+3\right)=2p-4-q.
\]

If $p> \frac{q+4}{2}$, then in view of (\ref{eq:2.4}) we have $\phi_{D_{p,q}^{0,1}}(P_3)> \phi_{B_{p,q}^{2,0}}(P_3)$. By Lemma 1.2(i),  $S_4(D_{p,q}^{0,1})>S_4(B_{p,q}^{2,0})$ holds. So we have $B_{p,q}^{2,0}\prec_sD_{p,q}^{0,1}$. So in this case $D_{p,q}^{0,1}$ is the third last tree, in the $S$-order, among $\mathscr{T}_n^{p,q}$. 

If $p= \frac{q+4}{2}$, then in view of (\ref{eq:2.4}) we have $\phi_{D_{p,q}^{0,1}}(P_3)=\phi_{B_{p,q}^{2,0}}(P_3)$. Hence, $S_4(D_{p,q}^{0,1})=S_4(B_{p,q}^{2,0})$ holds by Lemma 1.2(i). In view of Lemma 1.2(ii), $S_5(D_{p,q}^{0,1})=S_5(B_{p,q}^{2,0})$ holds. By direct computing, we have
\begin{eqnarray*}
  \phi_{D_{p,q}^{0,1}}(P_4)-\phi_{B_{p,q}^{2,0}}(P_4)&=&[(p-1)\times 1+(q-2)(p-1)]-[(p-3)(q-1)+2\times (q-1)]=0, \\
  \phi_{D_{p,q}^{0,1}}(K_{1,3})-\phi_{B_{p,q}^{2,0}}(K_{1,3})&=&\left({p \choose 3}+{{q-1} \choose 3}\right)-\left({{p-2}\choose 3}+{q\choose 3}+1\right)
                                                             =\frac{-(q-3)^{2}+1}{4}<0.
\end{eqnarray*}
The last inequality follows by $q>p\geqslant  4$. In view of Lemma 1.2(iii), we have $S_6(D_{p,q}^{0,1})-S_6(B_{p,q}^{2,0})=3[-(q-3)^{2}+1]<0$, i.e., $D_{p,q}^{0,1}\prec_s B_{p,q}^{2,0}$. That is to say, $B_{p,q}^{2,0}$ is the third last tree, in the $S$-order, among $\mathscr{T}_n^{p,q}$ for $p= \frac{q+4}{2}$.

If $p< \frac{q+4}{2}$, then in view of (\ref{eq:2.4}) we have $\phi_{D_{p,q}^{0,1}}(P_3)< \phi_{B_{p,q}^{2,0}}(P_3)$. By Lemma 1.2(i),  $S_4(D_{p,q}^{0,1})<S_4(B_{p,q}^{2,0})$ holds. So we have $D_{p,q}^{0,1}\prec_sB_{p,q}^{2,0}$. So $B_{p,q}^{2,0}$ is the third last tree, in the $S$-order, among $\mathscr{T}_n^{p,q}$ for $p< \frac{q+4}{2}$.

This completes the proof of (ii).
\end{proof}
\begin{figure}[h!]
\begin{center}
  \psfrag{1}{$k$}\psfrag{9}{$k'$}\psfrag{2}{$F_{p,q}^k$}\psfrag{3}{$L_{p,q}^{k'}$}\psfrag{4}{$M_{p,q}^{k'}$}\psfrag{5}{$N_{p,q}^k$}
  \includegraphics[width=140mm]{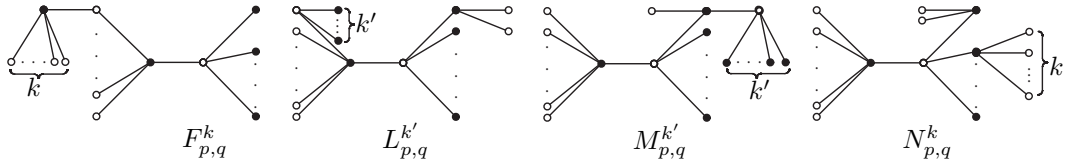}\\
  \caption{Trees $F_{p,q}^k, L_{p,q}^{k'}, M_{p,q}^{k'}$ and $N_{p,q}^k$ each of which contains $p$ white points and $q$ black points.}
\end{center}
\end{figure}

For convenience, let
$F_{p,q}^k,\, L_{p,q}^{k'},\,M_{p,q}^{k'}$ and $N_{p,q}^k\, (1\leqslant k\leqslant p-2, 1\leqslant k'\leqslant q-2)$ be the trees as depicted in Fig. 7, it is easy to see that $F_{p,q}^k,\, L_{p,q}^{k'},\,M_{p,q}^{k'},\, N_{p,q}^k$ are in $\mathscr{T}_n^{p,q}$.
\begin{thm}\label{lem2.6}
Given positive integers $p,q$ with $4\leqslant p< q$ and $p+q= n$.
\begin{wst}
  \item[{\rm (i)}] If $p> \frac{q+4}{2}$, then for any  $T\in \mathscr{T}_n^{p,q}\setminus \{B_{p,q}^{0,0},\,B_{p,q}^{0,1},\,D_{p,q}^{0,1}\}$, we have $T\preceq_s B_{p,q}^{2,0}$ with equality if and only if $T\cong B_{p,q}^{2,0}.$
  \item[{\rm (ii)}] If $p= \frac{q+4}{2}$, then for any  $T\in \mathscr{T}_n^{p,q}\setminus \{B_{p,q}^{0,0},\,B_{p,q}^{0,1},\,B_{p,q}^{2,0}\}$, we have $T\preceq_s D_{p,q}^{0,1}$ with equality if and only if $T\cong D_{p,q}^{0,1}.$
  \item[{\rm (iii)}] If $p< \frac{q+4}{2}$, then for any  $T\in \mathscr{T}_n^{p,q}\setminus \{B_{p,q}^{0,0},\,B_{p,q}^{0,1},\,B_{p,q}^{2,0}\}$, we have $T\preceq_s B_{p,q}^{0,2}$ with equality if and only if $T\cong B_{p,q}^{0,2}.$
\end{wst}
\end{thm}
\begin{proof}
For any $T \in \mathscr{T}_n^{p,q}$ such that $T\not\cong  B_{p,q}^{0,0},\, B_{p,q}^{0,1},\, D_{p,q}^{0,1},\, B_{p,q}^{2,0}$, by a similar discussion as above, $T$ can be transformed into $B_{p,q}^{0,0}$ (resp. $B_{p,q}^{0,1},\, D_{p,q}^{0,1},\, B_{p,q}^{2,0}$) by carrying Operations I and I\!I repeatedly. Let $\mathscr{C}_1$ (resp. $\mathscr{D}_1$) denote the set of all trees in $\mathscr{T}_n^{p,q}$ which can
be transformed into $D_{p,q}^{0,1}$ (resp. $B_{p,q}^{2,0}$) by carrying Operation I once, and let $\mathscr{C}_2$ (resp. $\mathscr{D}_2$) denote the set of all trees in $\mathscr{T}_n^{p,q}$ which can be transformed into $D_{p,q}^{0,1}$ (resp. $B_{p,q}^{2,0}$) by carrying Operation I\!I once.

(i)\ $p> \frac{q+4}{2}$. The last tree, in the $S$-order, among $\mathscr{T}_n^{p,q}\setminus \{B_{p,q}^{0,0},\,B_{p,q}^{0,1},\,D_{p,q}^{0,1}\}$ must be in
$\mathscr{A}_2\cup \mathscr{B}_1\cup \mathscr{B}_2\cup \mathscr{C}_1\cup \mathscr{C}_2$, where $\mathscr{A}_2$ is defined in the proof of Theorem 2.2, $\mathscr{B}_1,\, \mathscr{B}_2$ are defined in the proof of Theorem 2.3, while $\mathscr{C}_1= \{D_{p,q}^{2,0},\,D_{p,q}^{0,2},\,C_{p,q}^1,\,F_{p,q}^1\},\, \mathscr{C}_2=\{C_{p,q}^k: 2\leqslant k\leqslant q-2\} \cup \{F_{p,q}^k: 2\leqslant k \leqslant p-2\}.$

In view of the proof of Claims 1 and 2 in the proof of Theorem 2.3, we know that the last tree, in the $S$-order, among $\mathscr{A}_2\cup \mathscr{B}_1\cup \mathscr{C}_1$ is $B_{p,q}^{2,0}$. In what follows we show that for any $T$ in $\mathscr{C}_1\cup \mathscr{C}_2$, we have $T \prec_{s} B_{p,q}^{2,0}$.

In fact, $C_{p,q}^1$ (resp. $F_{p,q}^1$) can be obtained from $C_{p,q}^k$ (resp. $F_{p,q}^k$) by using {Operation I} $(k-1)$ times, where $k\geqslant  2$. By Lemma 1.4, we have $C_{p,q}^k\prec_s C_{p,q}^1$ and $F_{p,q}^k\prec_s F_{p,q}^1$ for $k\geqslant  2$. By the proof of Theorem 2.3, we know that $C_{p,q}^1\prec_s B_{p,q}^{2,0}$ and $D_{p,q}^{0,2}\prec_s D_{p,q}^{2,0}\prec_s B_{p,q}^{2,0}$. By Lemma 1.1, we have $S_i(B_{p,q}^{2,0})=S_i(F_{p,q}^1)$ for $i=0,1,2,3$. In view of Lemma 1.2(i), it is routine to check that $\phi_{B_{p,q}^{2,0}}(P_2)=\phi_{F_{p,q}^1}(P_2)=n-1,\, \phi_{B_{p,q}^{2,0}}(C_4)=\phi_{F_{p,q}^1}(C_4)=0.$ By Lemma 1.3,
$$
\phi_{B_{p,q}^{2,0}}(P_3)-\phi_{F_{p,q}^1}(P_3)={p-2\choose 2}+{q\choose 2}+{3\choose 2}-\left({p-1\choose 2}+{q-1\choose 2}+2\right)=q-p+2>0.
$$
Hence, $S_4(B_{p,q}^{2,0})-S_4(F_{p,q}^1)=4(q-p+2)>0,$ i.e., $F_{p,q}^1\prec_{s}B_{p,q}^{2,0}$. This completes the proof of (i).

\medskip
In what follows, we consider $p\leqslant \frac{q+4}{2}$.  By Lemmas 1,4, 1.5 and Theorem 2.3(ii), the last tree, in the $S$-order, among $\mathscr{T}_n^{p,q}\setminus \{B_{p,q}^{0,0},\,B_{p,q}^{0,1},\,B_{p,q}^{2,0}\}$ must be in
$\mathscr{A}_1\cup \mathscr{A}_2\cup \mathscr{B}_1\cup \mathscr{B}_2\cup \mathscr{D}_1\cup \mathscr{D}_2\setminus\{B_{p,q}^{0,1},\,B_{p,q}^{2,0}\}$, where $\mathscr{A}_1,\, \mathscr{A}_2$ are defined in the proof of Theorem 2.2, $\mathscr{B}_1,\, \mathscr{B}_2$ are defined in the proof of Theorem 2.3, while  $\mathscr{D}_2=\{L_{p,q}^k: 2\leqslant k\leqslant q-2\} \cup \{M_{p,q}^k: 2\leqslant k \leqslant q-2\}\cup \{N_{p,q}^k: 2\leqslant k \leqslant p-4\},$ $\mathscr{D}_1= \{B_{p,q}^{2,1},\,B_{p,q}^{0,2},\,L_{p,q}^1,\,M_{p,q}^1\}$ if $4\leqslant p<7$ and $\mathscr{D}_1= \{B_{p,q}^{3,0},\, B_{p,q}^{2,1},\,B_{p,q}^{0,2},\,L_{p,q}^1,\,M_{p,q}^1\}$ if $p\geqslant  7.$ 
\medskip

(ii)\ \  $p =\frac{q+4}{2}.$ In this case, we consider the following two subcases according to $\mathscr{D}_1.$

\medskip
{\bf Case 1.} $\mathscr{D}_1= \{B_{p,q}^{2,1},\,B_{p,q}^{0,2},\,L_{p,q}^1,\,M_{p,q}^1\}$ with $4\leqslant p < 7$.\medskip

First we determine the last tree, in the $S$-order, among $\mathscr{D}_1\cup \mathscr{D}_2.$ It is easy to see (based on Lemma 1.4), we have $B_{p,q}^{2,1}\prec_{s}B_{p,q}^{0,2}$.
Note that $L_{p,q}^1$ (resp. $M_{p,q}^1,\, B_{p,q}^{2,1}$) can be obtained from $L_{p,q}^k$ (resp. $M_{p,q}^k,\, N_{p,q}^k$) by using {Operation I}\ $(k-1)$ times, hence by Lemma 1.4 we have $L_{p,q}^k\prec_s L_{p,q}^1,\, M_{p,q}^k\prec_s M_{p,q}^1$ and $N_{p,q}^k\prec_s B_{p,q}^{2,1}$ for $k\geqslant  2$.

By Lemma 1.1, $S_i(L_{p,q}^1)=S_i(M_{p,q}^1)=S_i(B_{p,q}^{0,2})$ holds for $i=0,1,2,3$. By Lemma 1.2(i), we have $\phi_{L_{p,q}^1}(P_2)=\phi_{M_{p,q}^1}(P_2)=\phi_{B_{p,q}^{0,2}}(P_2)=n-1,\, \phi_{L_{p,q}^1}(C_{4})=\phi_{M_{p,q}^1}(C_4)=\phi_{B_{p,q}^{0,2}}(C_{4})=0.$
Note that $L_{p,q}^1$ and $M_{p,q}^1$ have the same degree sequence, by Lemma 1.3 $\phi_{L_{p,q}^1}(P_{3})=\phi_{M_{p,q}^1}(P_{3})$. Hence,
\begin{eqnarray*}
  \phi_{L_{p,q}^1}(P_{3})-\phi_{B_{p,q}^{0,2}}(P_{3}) &=& \phi_{M_{p,q}^1}(P_{3})-\phi_{B_{p,q}^{0,2}}(P_{3}) \\
   &=& \left({p-2 \choose 2}+{q-1\choose 2}+3+1\right)-\left({{p-2}\choose 2}+{q\choose 2}+2\right)\\
                                       &=&3-q<0.
\end{eqnarray*}
The last inequality follows by $q > p\geqslant  4$. By Lemma 1.2(i), we have $S_{4}(L_{p,q}^1)-S_{4}(B_{p,q}^{0,2})=S_{4}(M_{p,q}^1)-S_{4}(B_{p,q}^{0,2})=4(\phi_{M_{p,q}^1}(P_{3})-\phi_{B_{p,q}^{0,2}}(P_{3}))<0$, i.e., $L_{p,q}^1\prec_{s}B_{p,q}^{0,2}$ and $M_{p,q}^1\prec_{s}B_{p,q}^{0,2}$. Hence, $B_{p,q}^{0,2}$ is the last tree, in the $S$-order, among $\mathscr{D}_1\cup \mathscr{D}_2.$

By the proof of Claim 2 in Theorem 2.3, we obtain that $B_{p,q}^{0,2}$ is the last graph, in the $S$-order, among $(\mathscr{B}_1\cup \mathscr{B}_2)\setminus\{B_{p,q}^{2,0}\}.$

Note that $p<7$, it is routine to check that $(\mathscr{A}_1\cup \mathscr{A}_2)\setminus\{B_{p,q}^{0,1},\,B_{p,q}^{2,0}\}=\{D_{p,q}^{2,0},\,D_{p,q}^{0,1}\}.$
By Lemma 1.4, we have $D_{p,q}^{2,0} \prec_s D_{p,q}^{0,1}.$ In order to complete the proof, it suffices to compare $B_{p,q}^{0,2}$ with $D_{p,q}^{0,1}.$

By Lemma 1.1, $S_i(D_{p,q}^{0,1})=S_i(B_{p,q}^{0,2})$ holds for $i=0,1,2,3$. It is routine to check that $\phi_{D_{p,q}^{0,1}}(P_2)=\phi_{B_{p,q}^{0,2}}(P_2)=n-1$ and $\phi_{D_{p,q}^{0,1}}(C_4)=\phi_{B_{p,q}^{0,2}}(C_4)=0$. By Lemma 1.3,
$$
   \phi_{D_{p,q}^{0,1}}(P_3)-\phi_{B_{p,q}^{0,2}}(P_3)=\left({p\choose 2}+{q-1\choose 2}+1\right)-\left({q\choose 2}+{p-2\choose 2}+2\right)=2p-q-3=1.
$$
In view of Lemma 1.2(i), we have $S_4(B_{p,q}^{0,2})<S_4(D_{p,q}^{0,1})$, i.e.,
\[
  B_{p,q}^{0,2}\prec_s D_{p,q}^{0,1}.
\]
That is to say, our result holds in this case.

\medskip
{\bf Case 2.} $\mathscr{D}_1= \{B_{p,q}^{3,0},\, B_{p,q}^{2,1},\,B_{p,q}^{0,2},\,L_{p,q}^1,\,M_{p,q}^1\}$ with $p\geqslant  7$. \medskip

First we determine the last tree, in the $S$-order, among $\mathscr{D}_1\cup \mathscr{D}_2.$ In fact, by a similar discussion in Case 1 of determining the last graph, in the $S$-order, among $\mathscr{D}_1\cup \mathscr{D}_2,$ we can obtain that in this case, the last graph, in the $S$-order, among $(\mathscr{D}_1\setminus\{B_{p,q}^{3,0}\})\cup \mathscr{D}_2$ is just $B_{p,q}^{0,2}$. Hence, it suffices to compare $B_{p,q}^{3,0}$ with $B_{p,q}^{0,2}.$

In fact, by Lemma 1.1 $S_i(B_{p,q}^{0,2})=S_i(B_{p,q}^{3,0})$ holds for $i=0,1,2,3$. It is routine to check that $\phi_{B_{p,q}^{3,0}}(P_2)=\phi_{B_{p,q}^{0,2}}(P_2)=n-1$ and $\phi_{B_{p,q}^{3,0}}(C_4)=\phi_{B_{p,q}^{0,2}}(C_4)=0.$ By Lemma 1.3 we have
$$
   \phi_{B_{p,q}^{0,2}}(P_3)-\phi_{B_{p,q}^{3,0}}(P_3)=\left({p-2\choose 2}+{q\choose 2}+2\right)-\left({p-3\choose 2}+{q\choose 2}+{4\choose 2}\right)=p-7\geqslant 0.
$$

If $p>7$, by Lemma 1.2(i) $S_4(B_{p,q}^{0,2})>S_4(B_{p,q}^{3,0})$, i.e., $B_{p,q}^{3,0}\prec_s B_{p,q}^{0,2}$.

If $p=7$, we have $\phi_{B_{p,q}^{0,2}}(P_3)=\phi_{B_{p,q}^{3,0}}(P_3).$ By direct computing, we have $\phi_{B_{p,q}^{0,2}}(P_4) = \phi_{B_{p,q}^{3,0}}(P_4)=(p-1)(q-1)$ and
$$
   \phi_{B_{p,q}^{0,2}}(K_{1,3})-\phi_{B_{p,q}^{3,0}}(K_{1,3})={{p-2} \choose 3}+{q \choose 3}-{{p-3}\choose 3}-{q\choose 3}
                                                             =\frac{1}{2}(p-3)(p-4)>0.
$$
By Lemma 1.2(iii), we have $S_6(B_{p,q}^{0,2})-S_6(B_{p,q}^{3,0})=6(p-3)(p-4)>0$, i.e.,
\[
   B_{p,q}^{3,0}\prec_s B_{p,q}^{0,2}.
\]
Hence, $B_{p,q}^{0,2}$ is the last graph, in the $S$-order, among $\mathscr{D}_1\cup \mathscr{D}_2$ in this case.

By the proof of Claim 2 in Theorem 2.3, we obtain that $B_{p,q}^{0,2}$ is the last graph, in the $S$-order, among $(\mathscr{B}_1\cup \mathscr{B}_2)\setminus\{B_{p,q}^{2,0}\}.$

Note that $p\geqslant  7$, it is routine to check that
\[
(\mathscr{A}_1\cup \mathscr{A}_2)\setminus\{B_{p,q}^{0,1},\,B_{p,q}^{2,0}\}=\{D_{p,q}^{0,1}\}\cup\left\{B_{p,q}^{k,0}: 3\leqslant k\leqslant \left\lfloor\frac{p-1}{2}\right\rfloor\right\}\cup \left\{D_{p,q}^{k,0}: 2\leqslant k\leqslant\left\lfloor\frac{q-1}{2}\right\rfloor\right\}.
\]
By Lemma 1.4, we have $D_{p,q}^{2,0} \prec_s D_{p,q}^{0,1}.$ In view of (2.1), (2.2) and (2.7), it suffices to compare $B_{p,q}^{3,0}$ with $D_{p,q}^{0,1}.$

In view of (2.6), we obtain that $B_{p,q}^{3,0}\prec_s B_{p,q}^{0,2}$. When $p\geqslant  7$, by a similar discussion as in the proof of (2.5), we can also show that
$B_{p,q}^{0,2}\prec_s D_{p,q}^{0,1}.$ Hence,  $B_{p,q}^{3,0}\prec_s D_{p,q}^{0,1}.$

Combining with the proof as above, we obtain that $D_{p,q}^{0,1}$ is the fourth last tree, in the $S$-order, among $\mathscr{T}_n^{p,q}$.
This completes the proof of (ii).\medskip

(iii) Let $p< \frac{q+4}{2}$. We proceed by considering the following two possible cases according to $\mathscr{D}_1$.

\medskip
{\bf Case 1.}
$\mathscr{D}_1= \{ B_{p,q}^{2,1},\,B_{p,q}^{0,2},\,L_{p,q}^1,\,M_{p,q}^1\}$ with $4\leqslant p<7$. \medskip

By a similar discussion as the proof of Case 1 in (ii), we know that $B_{p,q}^{0,2}$ is the last tree, in the $S$-order, among $(\mathscr{A}_2\cup\mathscr{B}_1\cup
\mathscr{B}_2\cup\mathscr{D}_1\cup \mathscr{D}_2)\setminus\{B_{p,q}^{2,0}\}.$
Note that $p<7$, it is routine to check that $\mathscr{A}_1\setminus\{B_{p,q}^{0,1}\}=\{D_{p,q}^{0,1}\}.$
In order to complete the proof, it suffices to compare $B_{p,q}^{0,2}$ with $D_{p,q}^{0,1}.$

By Lemma 1.1, $S_i(D_{p,q}^{0,1})=S_i(B_{p,q}^{0,2})$ holds for $i=0,1,2,3$. It is routine to check that $\phi_{D_{p,q}^{0,1}}(P_2)=\phi_{B_{p,q}^{0,2}}(P_2)=n-1$ and
$\phi_{D_{p,q}^{0,1}}(C_4)=\phi_{B_{p,q}^{0,2}}(C_4)=0$. By Lemma 1.3,
$$
  \phi_{D_{p,q}^{0,1}}(P_3)-\phi_{B_{p,q}^{0,2}}(P_3)=\left({p\choose 2}+{q-1\choose 2}+1\right)-\left({q\choose 2}+{p-2\choose 2}+2\right)=2p-q-3.
$$

If $p<\frac{q+3}{2}$, by Lemma 1.2(i), we have $S_4(D_{p,q}^{0,1})<S_4(B_{p,q}^{0,2})$, i.e., $D_{p,q}^{0,1}\prec_s B_{p,q}^{0,2}.$

If $p=\frac{q+3}{2}$, we have $S_{4}(D_{p,q}^{0,1})=S_{4}(B_{p,q}^{0,2})$. By Lemma 1.2(ii), $S_5(D_{p,q}^{0,1})=S_5(B_{p,q}^{0,2})$.
By direct computing, we have $\phi_{D_{p,q}^{0,1}}(P_{4})=\phi_{B_{p,q}^{0,2}}(P_{4})=(p-1)(q-1)$ and
$$
  \phi_{D_{p,q}^{0,1}}(K_{1,3})-\phi_{B_{p,q}^{0,2}}(K_{1,3})={p\choose 3}+{q-1\choose 3}-{p-2\choose 3}-{q\choose 3}=\frac{-(q-2)^{2}+1}{4}<0.
$$
Hence, by Lemma 1.2(iii), we have $S_{6}(D_{p,q}^{0,1})-S_{6}(B_{p,q}^{0,2})=3[-(q-2)^{2}+1]<0$, i.e., $D_{p,q}^{0,1}\prec_{s}B_{p,q}^{0,2}$.
So in this case, $B_{p,q}^{0,2}$ is the fourth last tree, in the $S$-order, among $\mathscr{T}_n^{p,q}$.

\medskip
{\bf Case 2.} $\mathscr{D}_1= \{B_{p,q}^{3,0},\, B_{p,q}^{2,1},\,B_{p,q}^{0,2},\,L_{p,q}^1,\,M_{p,q}^1\}$  with $p\geqslant  7$. \medskip

By a similar discussion as the proof of Case 2 in (ii), we know that $B_{p,q}^{0,2}$ is the last tree, in the $S$-order, among $(\mathscr{A}_2\cup\mathscr{B}_1\cup
\mathscr{B}_2\cup\mathscr{D}_1\cup \mathscr{D}_2)\setminus\{B_{p,q}^{2,0}\}.$ It is routine to check that
$\mathscr{A}_1\setminus\{B_{p,q}^{0,1}\}=\{D_{p,q}^{0,1}\}.$
In order to complete the proof, it suffices to compare $B_{p,q}^{0,2}$ with $D_{p,q}^{0,1}.$ By a similar discussion as the proof of Case 1 in (iii), we have
$D_{p,q}^{0,1}\prec_{s}B_{p,q}^{0,2}$. Hence, in this case $B_{p,q}^{0,2}$ is the fourth last tree in the $S$-order, among $\mathscr{T}_n^{p,q}$.
This completes the proof of (iii).
\end{proof}

\begin{thm}\label{cor2.7}
If $4\leqslant p=q$, then for any $T\in \mathscr{T}_n^{p,q}\setminus \{B_{p,q}^{0,0},\, B_{p,q}^{0,1},\,B_{p,q}^{2,0}\}$, we have $T\preceq_s B_{p,q}^{0,2}$ with equality if and only if $T \cong B_{p,q}^{0,2}$.
\end{thm}
\begin{proof}
Up to isomorphism, for any $T \in \mathscr{T}_n^{p,q}$ such that $T\not\cong  B_{p,q}^{0,0},\, B_{p,q}^{0,1},\, B_{p,q}^{2,0}$, by a similar discussion as above, $T$ can be transformed into $B_{p,q}^{0,0}$ (resp. $B_{p,q}^{0,1},\, B_{p,q}^{2,0}$) by carrying the Operations I and I\!I repeatedly. By Lemmas 1.4 and 1.5, the last tree, in the $S$-order, among $\mathscr{T}_n^{p,q}\setminus \{B_{p,q}^{0,0},\,B_{p,q}^{0,1},\,B_{p,q}^{2,0}\}$ must be in
$\mathscr{A}_1\cup \mathscr{A}_2\cup \mathscr{B}_1\cup \mathscr{B}_2\cup \mathscr{D}_1\cup \mathscr{D}_2\setminus\{B_{p,q}^{0,1},\,B_{p,q}^{2,0}\}$, where $\mathscr{A}_1,\, \mathscr{A}_2$ (resp. $\mathscr{B}_1,\, \mathscr{B}_2$) are defined in the proof of Theorem 2.2 (resp. Theorem 2.3), and $\mathscr{D}_1,\, \mathscr{D}_2$ are defined in the proof of Theorem 2.4. We proceed by considering the following two possible cases.

\medskip
{\bf Case 1.} $\mathscr{D}_1= \{ B_{p,q}^{2,1},\,B_{p,q}^{0,2},\,L_{p,q}^1,\,M_{p,q}^1\}$ with $4\leqslant p <7$. \medskip

By a similar discussion as the proof of Case 1 in Theorem 2.4(ii), we obtain that $B_{p,q}^{0,2}$ is the last tree, in the $S$-order, among $\mathscr{D}_1\cup \mathscr{D}_2.$ By the proof of Claim 2 in Theorem 2.3, we obtain that $B_{p,q}^{0,2}$ is the last graph, in the $S$-order, among $(\mathscr{B}_1\cup \mathscr{B}_2)\setminus\{B_{p,q}^{2,0}\}.$ It is routine to check that $(\mathscr{A}_1\cup \mathscr{A}_2)\setminus\{B_{p,q}^{0,1},\,B_{p,q}^{2,0}\}=\emptyset$ in this case. Hence, $B_{p,q}^{0,2}$ is the last tree, in the $S$-order, among $\mathscr{T}_n^{p,q}\setminus \{B_{p,q}^{0,0},\, B_{p,q}^{0,1},\,B_{p,q}^{2,0}\}.$

\medskip
{\bf Case 2.} $\mathscr{D}_1= \{B_{p,q}^{3,0},\, B_{p,q}^{2,1},\,B_{p,q}^{0,2},\,L_{p,q}^1,\,M_{p,q}^1\}$ with $p \geqslant  7$. \medskip

By a similar discussion as the proof of Case 2 in Theorem 2.4(ii), we obtain that $B_{p,q}^{0,2}$ is the last tree, in the $S$-order, among $(\mathscr{D}_1\cup \mathscr{D}_2\cup \mathscr{A}_1\cup \mathscr{A}_2)\setminus\{B_{p,q}^{0,1},\,B_{p,q}^{2,0}\}$. By the proof of Claim 2 in Theorem 2.3, we obtain that $B_{p,q}^{0,2}$ is the last graph, in the $S$-order, among $(\mathscr{B}_1\cup \mathscr{B}_2)\setminus\{B_{p,q}^{2,0}\}.$ Hence, $B_{p,q}^{0,2}$ is the last tree, in the $S$-order, among $\mathscr{T}_n^{p,q}\setminus \{B_{p,q}^{0,0},\, B_{p,q}^{0,1},\,B_{p,q}^{2,0}\}.$

This completes the proof.
\end{proof}
\section{\normalsize Conclusion and remarks}
Summarize the results in Section 2, we can obtain the last four graphs in the $S$-order of the set of $n$-vertex trees with a $(p,q)$-bipartition.

Combining with Theorems 2.1, 2.2, 2.3(ii) and 2.4, we have
\begin{thm}\label{lem2.8}
Given positive integers $p,q$ with $4\leqslant p< q$ and $p+q= n$.
\begin{wst}
  \item[{\rm (i)}]If $p> \frac{q+4}{2},$ the last four trees, in the $S$-order, among $\mathscr{T}_n^{p, q}$ are as follows:
$
   B_{p,q}^{2,0}, D_{p,q}^{0,1}, B_{p,q}^{0,1}, B_{p,q}^{0,0}.
$
  \item[{\rm (ii)}]If $p= \frac{q+4}{2},$ the last four trees, in the $S$-order, among $\mathscr{T}_n^{p, q}$ are as follows:
$
   D_{p,q}^{0,1}, B_{p,q}^{2,0}, B_{p,q}^{0,1}, B_{p,q}^{0,0}.
$
  \item[{\rm (iii)}]If $p< \frac{q+4}{2},$ the last four trees, in the $S$-order, among $\mathscr{T}_n^{p, q}$ are as follows:
$
  B_{p,q}^{0,2}, B_{p,q}^{2,0}, B_{p,q}^{0,1}, B_{p,q}^{0,0}.
$
\end{wst}
\end{thm}

Combining with Theorems 2.1, 2.2, 2.3(i) and 2.5, we have
\begin{thm}\label{lem2.9}
If $4\leqslant p=q$, the last four trees, in the $S$-order, among the set $\mathscr{T}_n^{p, q}$ are as follows:
$
   B_{p,q}^{0,2},\,B_{p,q}^{2,0},\,B_{p,q}^{0,1},\,B_{p,q}^{0,0}.
$
\end{thm}

In this paper, we determine the the last four graphs, in the $S$-order, of the set of $n$-vertex trees with a $(p,q)$-bipartition. It is natural to consider
the following research problem: How can we determine the first $k$ graphs, in the $S$-order, of the set of $n$-vertex trees with a $(p,q)$-bipartition? It seems difficult but interesting.

\end{document}